\documentclass[11pt]{article} 
 
\usepackage{latexsym} 
\usepackage{amsfonts} 
\usepackage{amsmath} 
\usepackage{amsthm} 
\usepackage{mathrsfs} 
\usepackage{dsfont} 
\usepackage{bbold} 
\usepackage[english]{babel} 
\usepackage{caption} 
\usepackage{epsfig} 
\usepackage{float} 
\usepackage{subfigure} 
\usepackage{psfrag} 
\usepackage{graphicx} 
\usepackage{epsfig} 
\usepackage{amsbsy} 
\usepackage{hyperref}
 
\DeclareMathOperator{\Val}{\matV}

\newtheorem{theorem}{Theorem} 
\newtheorem*{prop*}{Theorem}

\newtheorem{lemma}[theorem]{Lemma}

\newtheorem{hyp}{Hypothesis}

\newcommand{\zerarcounters}{\setcounter{equation}{0}\setcounter{theorem}{0}} 
 
\newcommand{\ZZZ}{\mathds{Z}} 
\newcommand{\CCC}{\mathds{C}} 
\newcommand{\NNN}{\mathds{N}} 
 
\newcommand{\RRR}{\mathds{R}} 
\newcommand{\TTT}{\mathds{T}}

\newcommand{\BB}{{\mathcal B}}

\newcommand{\calG}{{\mathcal G}}

\newcommand{\TT}{{\mathcal T}}

\newcommand{\gotB}{{\mathfrak B}}

\newcommand{\gotT}{{\mathfrak T}}

\newcommand{\matV}{{\mathscr V}}

\newcommand{\ol}{\overline} 
 
\newcommand{\Fullbox}{{\rule{2.0mm}{2.0mm}}} 
 
\newcommand{\EP}{\hfill\Fullbox\vspace{0.2cm}} 
\newcommand{\prova}{\noindent{\it Proof. }} 
\newcommand{\io}{\infty} 
\newcommand{\e}{\varepsilon} 
\newcommand{\al}{\alpha}

\newcommand{\x}{\xi}

\newcommand{\g}{\gamma}

\newcommand{\oo}{\boldsymbol{\omega}} 
\newcommand{\mm}{\boldsymbol{\mu}}

\newcommand{\nn}{\boldsymbol{\nu}} 

\newcommand{\pps}{\boldsymbol{\psi}} 
\newcommand{\vzero}{\boldsymbol{0}}

\newcommand{\ii}{{\rm i}}

\def\ins#1#2#3{\vbox to0pt{\kern-#2 \hbox{\kern#1 #3}\vss}\nointerlineskip} 
 
\begin{document}
 
%%%%%%%%%%%%%%%%%%%%%%%%%%%%%%%%%%%%%%%%%%%%%%%%%
\title{\bf Domains of analyticity for response solutions\\
in strongly dissipative forced systems}
%%%%%%%%%%%%%%%%%%%%%%%%%%%%%%%%%%%%%%%%%%%%%%%%%
 
\author 
{\bf Livia Corsi$^{1}$, Roberto Feola$^{2}$ and Guido Gentile$^{3}$
\vspace{2mm} 
\\ \small
$^{1}$ Dipartimento di Matematica, Universit\`a di
Napoli ``Federico II'', Napoli, I-80126, Italy
\\ \small 
$^{2}$ Dipartimento di Matematica, Universit\`a di Roma ``La Sapienza",
Roma, I-00185, Italy
\\ \small 
$^{3}$ Dipartimento di Matematica, Universit\`a di Roma Tre, Roma,
I-00146, Italy
\\ \small 
E-mail:  livia.corsi@unina.it, feola@mat.uniroma1.it,
gentile@mat.uniroma3.it}

\date{} 
 
\maketitle 
 
%%%%%%%%%%%%%%%%%%%%%%%%%%%%%%%%%%%%%%%%%%%%%%%%%
\begin{abstract} 
We study the ordinary differential equation $\e\ddot x + \dot x + \e g(x) = \e f(\oo t)$,
where $g$ and $f$ are real-analytic functions, with $f$ quasi-periodic in $t$
with frequency vector $\oo$.  If $c_{0} \in\RRR$ is such that $g(c_0)$
equals the average of $f$ and $g'(c_0)\neq0$, under very mild assumptions
on $\oo$ there exists a quasi-periodic solution close to $c_0$
with frequency vector $\oo$.
We show that such a solution depends analytically on $\e$ in
a domain of the complex plane tangent more than quadratically
to the imaginary axis at the origin.
\end{abstract} 
%%%%%%%%%%%%%%%%%%%%%%%%%%%%%%%%%%%%%%%%%%%%%%%%%
  
%%%%%%%%%%%%%%%%%%%%%%%%%%%%%%%%%%%%%%%%%%%%%%%%%
%%%%%%%%%%%%%%%%%%%%%%%%%%%%%%%%%%%%%%%%%%%%%%%%%
\zerarcounters 
\section{Introduction} 
\label{sec:1} 
%%%%%%%%%%%%%%%%%%%%%%%%%%%%%%%%%%%%%%%%%%%%%%%%%
%%%%%%%%%%%%%%%%%%%%%%%%%%%%%%%%%%%%%%%%%%%%%%%%%

Consider the ordinary differential equation in $\RRR$
\begin{equation}\label{eq:1.1}
\e \ddot x + \dot x + \e \, g(x) = \e \, f(\oo t) ,
\end{equation}
where $\e\in\RRR$ is small and $\oo \in \RRR^{d}$, with $d\in\NNN$,
is assumed (without loss of generality) to have
rationally independent components, i.e.
$\oo\cdot\nn\neq0$ $\forall\nn\in\ZZZ^{d}_{*}:=\ZZZ^{d}\setminus\{\vzero\}$.
For $\e>0$ the equation describes a one-dimensional system
with mechanical force $g$, subject to a quasi-periodic forcing $f$
with frequency vector $\oo$ and in the presence of strong dissipation.
We refer to \cite{GBD1} for some physical background. 
A quasi-periodic solution to (\ref{eq:1.1}) with the
same frequency vector $\oo$ as the forcing will be called a \textit{response solution}.

%%%%%%%%%%%%%%%%%%%%%%%%%%%%%%%%%%%%%%%%%%%%%%%%%
\begin{hyp} \label{hyp:1}
The functions $g\!:\RRR\to\RRR$ and $f\!:\TTT^{d}\to\RRR$ are real-analytic. 
There is $c_{0}\in\RRR$ such that $g(c_{0})=f_{\vzero}$,
where $f_{\vzero}$ is the average of $f$ on $\TTT^{d}$, and $a:=g'(c_{0})\neq0$.
\end{hyp}
%%%%%%%%%%%%%%%%%%%%%%%%%%%%%%%%%%%%%%%%%%%%%%%%%

In other words we assume that $c_{0}$ is a simple zero of the function $g(x)-f_{\vzero}$.
Denote by $\Sigma_{\xi}:=\{ \pps=(\psi_{1},\ldots,\psi_{d})\in(\CCC/2\pi \ZZZ)^{d}  :
|{\rm Im}\,\psi_{k}| \le \xi \hbox{ for }k=1,\ldots,d\}$, with $\xi>0$,
the strip where $f$ is analytic. By the analyticity assumptions one can write
\begin{equation} \nonumber
f(\pps) = \sum_{\nn\in\ZZZ^{d}} {\rm e}^{\ii\nn\cdot\pps} f_{\nn} , \qquad 
g(x) =\sum_{p=0}^{\io} a_{p} \left( x - c_0 \right)^{p} ,
\end{equation}
where
\begin{equation} \nonumber
|f_{\nn}| \le \Phi \, {\rm e}^{-\x |\nn|} , \qquad
a_{p} := \frac{1}{p!} \frac{{\rm d}^{p}g}{{\rm d}x^{p}}(c_{0}) , \qquad
|a_{p}| \le \Gamma \, \rho^{\,p} ,
\end{equation}
for suitable constants $\Phi$, $\Gamma$ and $\rho$.
Set $N(f)=N$ if $f$ is a trigonometric polynomial of degree $N$ and
$N(f)=\io$ otherwise, and define
\begin{eqnarray} 
& & \beta_{n}(\oo) := \min \big\{ |\oo\cdot\nn| : 0<|\nn| \le 2^{n} ,|\nn|\le N(f) \big\} ,
\qquad \e_{n}(\oo) := \frac{1}{2^{n}} \log \frac{1}{\beta_{n}(\oo)} ,
\nonumber \\
& & \al_{n}(\oo) := \min \big\{ |\oo\cdot\nn| : 0<|\nn| \le 2^{n}  \big\} ,
\qquad \gotB(\oo) := \sum_{n=0}^{\io} \frac{1}{2^{n}} \log \frac{1}{\al_{n}(\oo)} . 
\nonumber
\end{eqnarray}
%

%%%%%%%%%%%%%%%%%%%%%%%%%%%%%%%%%%%%%%%%%%%%%%%%%
\begin{hyp} \label{hyp:2}
$\displaystyle{\lim_{n\to\io}\e_{n}(\oo)=0}$.
\end{hyp}
%%%%%%%%%%%%%%%%%%%%%%%%%%%%%%%%%%%%%%%%%%%%%%%%%

In particular no assumption at all is required on $\oo$ if $f$ is a trigonometric polynomial,
since $\beta_{n}(\oo)$ is definitively constant in that case.

Before stating our results we need some more notations. We define the sets
$C_{R} :=  \{\e\in\CCC : |{\rm Re}\,\e^{-1}| > (2R)^{-1}\}$ and $\Omega_{R,B} 
:= \{ \e\in\CCC : |{\rm Re}\,\e| \ge B \, ({\rm Im}\,\e)^{2} \hbox{ and } 0<|\e|<2R\}$.
$C_{R}$ consists of two disks with radius $R$ and centers $(R,0)$ and $(-R,0)$, 
while $\Omega_{R,B}$ is the intersection of 
the disk of center $(0,0)$ and radius $2R$ with two parabolas
with vertex at the origin: all such sets are tangent at the origin to the imaginary axis.
Note that the smaller $B$, the more flattened are the parabolas. If $2RB<1$ one has
$C_{R}\subset \Omega_{R,B}$.

The following result has been proved in \cite{CCD}.

%%%%%%%%%%%%%%%%%%%%%%%%%%%%%%%%%%%%%%%%%%%%%%%%%
\begin{theorem} \label{thm:1}
Assume Hypotheses \ref{hyp:1} and \ref{hyp:2} for the system (\ref{eq:1.1}) and
denote by $\Sigma_{\xi}$ the strip of analyticity of $f$.
Then there exist $\e_{0}>0$ and $B_{0}>0$ such that for all $B>B_{0}$ there is
a response solution $x(t)=c_{0}+u(\oo t,\e)$ to (\ref{eq:1.1}), with
$u(\pps,\e)=O(\e)$ analytic in $\pps\in\Sigma_{\xi'}$ and $\e\in \Omega_{\e_{0},B}$,
for some $\xi'<\xi$.
\end{theorem}
%%%%%%%%%%%%%%%%%%%%%%%%%%%%%%%%%%%%%%%%%%%%%%%%%

In the theorem above $\e_0$ has to be small, while $B_{0}$ must be large enough.
However, for $B$ as close as wished to $B_{0}$ one can take $\ol{\e}<\e_{0}$
small enough for the condition $\ol{\e}B<1$ to be satisfied, so as to obtain
that $C_{\ol{\e}/2}$ is contained inside the analyticity domain. In this respect
Theorem \ref{thm:1} extends previous results in the literature \cite{GBD1,GBD2},
where analyticity in a pair of disks was obtained under stronger conditions
on $\oo$, such as the standard Diophantine condition 
\begin{equation} \label{eq:1.2}
\left| \oo\cdot\nn \right| > \frac{\g}{|\nn|^{\tau}} \qquad \forall \nn \in \ZZZ^{d}_{*} ,
\end{equation}
or the Bryuno condition $\BB(\oo)<\io$ 
If either $d=1$ or $d=2$ and $\oo$ satisfies the standard
Diophantine condition (\ref{eq:1.2}) with $\tau=1$,
the response solution is Borel-summable.

In the present letter we remove in Theorem \ref{thm:1}
the condition for $B$ to be large, by proving the following result.

%%%%%%%%%%%%%%%%%%%%%%%%%%%%%%%%%%%%%%%%%%%%%%%%%
\begin{theorem} \label{thm:2}
Assume Hypotheses \ref{hyp:1} and \ref{hyp:2} for the system (\ref{eq:1.1}) and
denote by $\Sigma_{\xi}$ the strip of analyticity of $f$.
Then for all $B>0$ there exists $\e_{0}>0$ such that there is
a response solution $x(t)=c_{0}+u(\oo t,\e)$ to (\ref{eq:1.1}), with
$u(\pps,\e)=O(\e)$ analytic in $\pps\in\Sigma_{\xi'}$ and $\e\in \Omega_{\e_{0},B}$,
for some $\xi'<\xi$.
The dependence of $\e_0$ on $B$ is of the form $\e_{0}=\e_{1}B^{\alpha}$,
for some $\alpha>0$ and $\e_1$ independent of $B$.
\end{theorem}
%%%%%%%%%%%%%%%%%%%%%%%%%%%%%%%%%%%%%%%%%%%%%%%%%

The proof of the theorem given in Section \ref{sec:3} yields the value $\alpha=8$:
such a value is non-optimal and could be improved by a more careful analysis.
Thanks to Theorem \ref{thm:2}  we can estimate the domain
of analyticity by the union of the domains $\Omega_{\e_{0},B}$,
with $\e_{0}=B^{\alpha}\e_{1}$, by letting $B$ varying in $(0,1]$.
This provides a domain that near the origin has boundary of the form
$|{\rm Re}\,\e| \approx \e_{1}^{-\beta}|{\rm Im}\,\e|^{2+\beta}$, where $\beta=1/\alpha$.

%%%%%%%%%%%%%%%%%%%%%%%%%%%%%%%%%%%%%%%%%%%%%%%%%
%%%%%%%%%%%%%%%%%%%%%%%%%%%%%%%%%%%%%%%%%%%%%%%%%
\zerarcounters 
\section{Tree representation} 
\label{sec:2} 
%%%%%%%%%%%%%%%%%%%%%%%%%%%%%%%%%%%%%%%%%%%%%%%%%
%%%%%%%%%%%%%%%%%%%%%%%%%%%%%%%%%%%%%%%%%%%%%%%%%

We can rewrite (\ref{eq:1.1}) as
\begin{equation} \label{eq:2.1}
\e \ddot x + \dot x + \e \, a \, (x - c_{0}) + \mu \e \sum_{p=2}^{\io}
a_{p} (x-c_{0})^{p}
= \mu \e \sum_{\nn\in\ZZZ^{d}_{*}} {\rm e}^{\ii\nn\cdot\pps} f_{\nn} ,
\end{equation}
where $a:=a_{1}$ and $\mu=1$. However, we can consider $\mu$ as a free parameter
and study (\ref{eq:2.1}) for $\e\in\CCC$ and $\mu\in\RRR$.
Then we look for a quasi-periodic solution to (\ref{eq:2.1}) of the form
\begin{equation} \label{eq:2.2}
x(t,\e,\mu) = c_{0} + u(\oo t, \e,\mu) , \qquad u(\pps,\e,\mu) =
\sum_{k=1}^{\io} \sum_{\nn\in\ZZZ^{d}} 
\mu^{k} {\rm e}^{\ii\nn\cdot\pps} u^{(k)}_{\nn}(\e) .
\end{equation}
By inserting (\ref{eq:2.2}) into (\ref{eq:2.1}) we obtain a recursive definition
for the coefficients $u^{(k)}_{\nn}(\e)$,
which admits a natural graphical representation in terms of trees.

A \textit{rooted tree} $\theta$ is a graph with no cycle,
such that all the lines are oriented toward a unique
point (\textit{root}) which has only one incident line (root line).
All the points in $\theta$ except the root are called \textit{nodes}.
The orientation of the lines in $\theta$ induces a partial ordering 
relation ($\preceq$) between the nodes. Given two nodes $v$ and $w$,
we shall write $w \prec v$ every time $v$ is along the path
(of lines) which connects $w$ to the root. 
We shall write $w\prec \ell$ if $w\preceq v$, where $v$ is the node which $\ell$ exits.
For any node $v$ denote by $p_{v}$ the number of lines entering $v$:
$v$ is called and \textit{end node} if $p_{v}=0$ and an \textit{internal node} if $p_{v}>0$.
We denote by $N(\theta)$ the set of nodes,
by $E(\theta)$ the set of end nodes, by $V(\theta)$ the set of internal nodes
and by $L(\theta)$ the set of lines; one has $N(\theta)=E(\theta) \amalg V(\theta)$.

We associate with each end node $v\in E(\theta)$
a \textit{mode} label $\nn_{v}\in\ZZZ^{d}_{*}$
and with each internal node an \textit{degree} label $d_{v}\in\{0,1\}$.
With each line $\ell\in L(\theta)$ we associate
a \textit{momentum} $\nn_{\ell} \in \ZZZ^{d}$.
We impose the following constraints on the labels:
\begin{enumerate}
\item $\nn_{\ell}=\sum_{w \in E_{\ell}(\theta)}  \nn_{w}$, 
where $E_{\ell}(\theta):=\{ w\in E(\theta) : w \prec \ell \}$;
\item $p_{v}\ge 2$ $\forall v\in V(\theta)$;
\item if $d_{v}=0$ then the line $\ell$ exiting $v$ has $\nn_{\ell}=\vzero$.
\end{enumerate}

We shall write $V(\theta)=V_{0}(\theta) \amalg V_{1}(\theta)$,
where $V_{0}(\theta):=\{ v\in V(\theta): d_{v}=0\}$.
For any discrete set $A$ we denote by $|A|$ its cardinality.
Define the \textit{degree} and the \textit{order} of $\theta$ as
$d(\theta):=|E(\theta)|+|V_{1}(\theta)|$ and $k(\theta):=|N(\theta)|$, respectively.

We call \textit{equivalent} two labelled rooted trees which can be transformed into
each other by continuously deforming the lines in such a way that
they do not cross each other. In the following we shall consider only
inequivalent labelled rooted trees, and we shall call them call trees \textit{tout court},
for simplicity.

We associate with each node $v\in N(\theta)$  a \textit{node factor} $F_{v}$ 
and with each line $\ell\in L(\theta)$ a \textit{propagator} $\calG_{\ell}$, such that
\begin{equation} \nonumber
F_{v} := \begin{cases}
- \e^{d_{v}}  \, a_{p_{v}} , & v \in V(\theta) , \\
\e \, f_{\nn_{v}} , & v \in E(\theta) ,
\end{cases}
\qquad\qquad
\calG_{\ell} := \begin{cases}
1/D(\e, \oo\cdot\nn_{\ell}) , & \nn_{\ell} \neq \vzero , \\
1/a , & \nn_{\ell}=\vzero ,
\end{cases}
\end{equation}
where $D(\e,s) := - \e s^{2} + \ii s + \e \, a$. Then, by defining
\begin{equation} \label{eq:2.3}
\Val(\theta,\e) := \Biggl( \prod_{v\in N(\theta)} 
F_{v} \Biggr) \Biggl( \prod_{\ell\in L(\theta)} \calG_{\ell} \Biggr)
\end{equation}
one has
\begin{equation} \label{eq:2.4}
u^{(k)}_{\nn}(\e) = \sum_{\theta\in \TT_{k,\nn}} \Val(\theta,\e) , \quad \nn
\in \ZZZ^{d}
% \neq 0 , \qquad
%u^{(k)}_{\vzero}(\e) = \sum_{\theta\in \TT_{k,\vzero}} \Val(\theta,\e) , 
\end{equation}
where $\TT_{k,\nn}$ is the set of trees of order $k$
and momentum $\nn$ associated with the root line.
Note that $u^{(1)}_{\vzero}=0$ and $u^{(2)}_{\nn}=0$ for all $\nn\in\ZZZ^{d}$.

%%%%%%%%%%%%%%%%%%%%%%%%%%%%%%%%%%%%%%%%%%%%%%%%%
%%%%%%%%%%%%%%%%%%%%%%%%%%%%%%%%%%%%%%%%%%%%%%%%%
\zerarcounters 
\section{Proof of Theorem \ref{thm:2}} 
\label{sec:3} 
%%%%%%%%%%%%%%%%%%%%%%%%%%%%%%%%%%%%%%%%%%%%%%%%%
%%%%%%%%%%%%%%%%%%%%%%%%%%%%%%%%%%%%%%%%%%%%%%%%%

We shall prove Theorem \ref{thm:2} in the case in which $N(f)=\io$. The case
of trigonometric polynomials is in fact easier and can be dealt with as shown in \cite{CFG}.

%%%%%%%%%%%%%%%%%%%%%%%%%%%%%%%%%%%%%%%%%%%%%%%%%
\begin{lemma} \label{lem:3.1}
Set $c_{0}=\min\{1/8, B/18,B/8|a|,|a|/8,|a|B/4,\sqrt{|a|}/2\}$.
There exists $\e_{1}>0$ such that one has
$|D(\e,s)| \ge c_{0} \max\{\min\{1,s^{2}\},|\e|^{2}\}$
for all $s\in\RRR$ and all $\e\in \Omega_{B,\e_{1}}$.
\end{lemma}
%%%%%%%%%%%%%%%%%%%%%%%%%%%%%%%%%%%%%%%%%%%%%%%%%

%%%%%%%%%%%%%%%%%%%%%%%%%%%%%%%%%%%%%%%%%%%%%%%%%
\prova
Write $\e=x+\ii y$, with $|x| \ge By^{2}$ and $x$ small enough.
By symmetry it is enough to study $y\ge0$.
One has $|D(\e,s)|^{2}=(s+ya-y s^{2})^{2} + x^{2}(a-s^{2})^{2}$.
If $y=0$ the bound is straightforward.
If $y>0$ denote by $s_{1}$ and $s_{2}$ the two roots of $s+ya-y s^{2}=0$:
one has $s_{1}=-ay+O(y^2)$ and $s_{2}=1/y+ay +O(y^2)$.
Let $\e_1$ be so small that $|s_{1}+ay|\le |a|y/2$,
$|s_{2}-1/y| \le 1/6y$ and $18|a|y^2 \le 1$ for $|\e|\le \e_{1}$.
The following inequalities are easily checked:
(1) if $|s|<2|a|y$, then $|x|\,|a-s^{2}| \ge |ax|/2 \ge |a|By^{2}/2 \ge Bs^{2}/8|a|$; 
(2) if $|s-s_{2}| < 1/2y$, then $|x|\,|a-s^{2}| \ge |x|s^{2}/2 \ge |x|/18y^{2} \ge B/18$;
(3) if $|s|\ge 2|a|y$ and $|s-s_{2}| \ge 1/2y$, then 
(3.1) $|s+ya-y s^{2}| \ge y|s-s_{1}|\,|s-s_{2}| \ge |a|y/4$,
(3.2) $|s+ya-y s^{2}| \ge |s-s_{1}|/2\ge |s|/8$,
(3.3) if either $a<0$ or $a>0$ and $|a-s^{2}| > |a|/2$
one has $|x|\,|a-s^{2}| > |ax|/2$,
while if $a>0$ and $|a-s^{2}| \le |a|/2$
one has $|s+ya-y s^2|\ge |s| - y|a-s^{2}| \ge \sqrt{a}/2$.
By collecting together all the bounds the assertion follows.
\EP
%%%%%%%%%%%%%%%%%%%%%%%%%%%%%%%%%%%%%%%%%%%%%%%%%

%%%%%%%%%%%%%%%%%%%%%%%%%%%%%%%%%%%%%%%%%%%%%%%%%
\begin{lemma} \label{lem:3.2}
For any tree $\theta$ one has $|E(\theta)|\ge |V(\theta)|+1$ and hence
$2|E(\theta)|\ge k(\theta)+1$.
\end{lemma}
%%%%%%%%%%%%%%%%%%%%%%%%%%%%%%%%%%%%%%%%%%%%%%%%%

%%%%%%%%%%%%%%%%%%%%%%%%%%%%%%%%%%%%%%%%%%%%%%%%%
\prova
By induction on the order $k(\theta)$.
\EP
%%%%%%%%%%%%%%%%%%%%%%%%%%%%%%%%%%%%%%%%%%%%%%%%%

For $v\in V_{1}(\theta)$ define
$E(\theta,v):=\{w\in E(\theta) : \hbox{ the line exiting $w$ enters $v$}\}$
and set $r_{v}:=|E(\theta,v)|$, $s_{v}:=p_{v}-r_{v}$,
$\mm_{v}:=\sum_{w\in E(\theta,v)} \nn_{w}$ and $\mu_{v}:=|\mm_{v}|$.
Define $V_{2}(\theta):=\{v\in V(\theta) : s_{v}=0\}$
and $V_{3}(\theta):=\{v\in V(\theta) : r_{v}=s_{v}=1\}$.
For $v\in V_{2}(\theta)$ call $\ell_{v}$ the line exiting $v$, and
for $v\in V_{3}(\theta)$ call $\ell_{v}$ the line exiting $v$ and $\ell_{v}'$
the line entering $v$ which does not exits an end node.
Define $\ol{V}_{2}(\theta) := \{v\in V_{2}(\theta) : \nn_{\ell_{v}}\neq\vzero\}$ and
$\ol{V}_{3}(\theta) := \{v\in V_{3}(\theta) : \nn_{\ell_{v}}\neq\vzero \hbox{ and }
\nn_{\ell_{v}'}\neq\vzero\}$, and set
$\ol{V}_{1}(\theta)=\ol{V}_{2}(\theta)\amalg \ol{V}_{3}(\theta)$.
By construction one has $\ol{V}_{1}(\theta) \subset V_{1}(\theta)$.

%%%%%%%%%%%%%%%%%%%%%%%%%%%%%%%%%%%%%%%%%%%%%%%%%
\begin{lemma} \label{lem:3.3}
\emph{
There exists $C_{0}>0$ such that
$C_{0}|\oo\cdot\nn| \ge {\rm e}^{-\xi|\nn|/16}$ $\forall \nn\in\ZZZ^{d}_{*}$.
}
\end{lemma}
%%%%%%%%%%%%%%%%%%%%%%%%%%%%%%%%%%%%%%%%%%%%%%%%%

%%%%%%%%%%%%%%%%%%%%%%%%%%%%%%%%%%%%%%%%%%%%%%%%%
\prova
It follows from Hypothesis \ref{hyp:2} by using that $\beta_{n}(\oo)=\alpha_{n}(\oo)$
if $N(f)=\io$.
\EP
%%%%%%%%%%%%%%%%%%%%%%%%%%%%%%%%%%%%%%%%%%%%%%%%%

%%%%%%%%%%%%%%%%%%%%%%%%%%%%%%%%%%%%%%%%%%%%%%%%%
\begin{lemma} \label{lem:3.4}
One has $C_{0}|\oo\cdot\nn_{\ell_{v}}| \ge {\rm e}^{-\xi \mu_{v}/16}$
for $v\in \ol{V}_{2}(\theta)$ and
$2 C_{0}\max\{|\oo\cdot\nn_{\ell_{v}}|,
|\oo\cdot\nn_{\ell_{v}'}|\} \ge {\rm e}^{-\xi \mu_{v}/16}$ for $v\in \ol{V}_{3}(\theta)$.
\end{lemma}
%%%%%%%%%%%%%%%%%%%%%%%%%%%%%%%%%%%%%%%%%%%%%%%%%

%%%%%%%%%%%%%%%%%%%%%%%%%%%%%%%%%%%%%%%%%%%%%%%%%
\prova
For $v\in \ol{V}_{2}(\theta)$ one has $\nn_{\ell_{v}}=\mm_{v}$,
so that one can use Lemma \ref{lem:3.3}.
For $v\in \ol{V}_{3}(\theta)$ one proceeds by contradiction.
Suppose that the assertion is false: this would imply
\begin{equation} \nonumber
{\rm e}^{-\xi \mu_{v}/16}  > C_{0} |\oo\cdot\nn_{\ell_{v}}| + 
C_{0} |\oo\cdot\nn_{\ell_{v}'}|
\ge C_{0} |\oo\cdot(\nn_{\ell_{v}}-\nn_{\ell_{v}'})|= C_{0}
| \oo\cdot \mm_{v} | \ge {\rm e}^{-\xi \mu_{v}/16} ,
\end{equation}
where we have used that $E(\theta,v)$ contains only one node $w$
and hence $\mm_{v} =\nn_{w} \neq \vzero$.
\EP
%%%%%%%%%%%%%%%%%%%%%%%%%%%%%%%%%%%%%%%%%%%%%%%%%

Define $L_{1}(\theta,v):=\{\ell_{v}\}$ for $v\in \ol{V}_{2}(\theta)$ and
$L_{1}(\theta,v):=\{\ell\in\{\ell_{v},\ell_{v}'\}:
2C_{0} |\oo\cdot\nn_{\ell}| \ge {\rm e}^{-\xi \mu_{v}/16}\}$ for $v\in \ol{V}_{3}(\theta)$.
Lemma \ref{lem:3.4} yields
$L_{1}(\theta,v)\neq\emptyset$ for all $v\in \ol{V}_{1}(\theta)$.
Set also $L_{1}(\theta):=\{ \ell \in L(\theta) : \exists v \in
\ol{V}_{1}(\theta)\hbox{ such that } \ell \in L_{1}(\theta,v)\}$,
$L_{\rm int}(\theta):=\{\ell\in L(\theta) : \ell \hbox{ exits
a node } v\in V_{1}(\theta)\}$ and $L_{0}(\theta) := L_{\rm int}(\theta)\setminus
L_{1}(\theta)$.

%%%%%%%%%%%%%%%%%%%%%%%%%%%%%%%%%%%%%%%%%%%%%%%%%
\begin{lemma} \label{lem:3.5}
For any tree $\theta$ one has $4 \, |L_{0}(\theta)| \le 3 |E(\theta)| - 4$.
\end{lemma}
%%%%%%%%%%%%%%%%%%%%%%%%%%%%%%%%%%%%%%%%%%%%%%%%%

%%%%%%%%%%%%%%%%%%%%%%%%%%%%%%%%%%%%%%%%%%%%%%%%%
\prova
By induction on $V(\theta)$.
If $|V(\theta)|=1$ then either $V(\theta)=V_{0}(\theta)$ or $V(\theta)=\ol{V}_{2}(\theta)$
and hence $|L_{0}(\theta)|=0$, so that the bound holds.
If $|V(\theta)|\ge 2$ the root line $\ell_{0}$ of $\theta$ exits a
node $v_{0}\in V(\theta)$ with $s_{v_{0}}+r_{v_{0}}\ge 2$
and $s_{v_{0}}\ge 1$. Call $\theta_{1},\ldots,\theta_{s_{v_{0}}}$
the trees whose respective root lines $\ell_{1},\ldots,\ell_{s_{v_{0}}}$
enter $v_{0}$: one has $|E(\theta)|=|E(\theta_{1})|+\ldots+|E(\theta_{s_{v_{0}}})|
+r_{v_{0}}$. If $\ell_{0}\notin L_{0}(\theta)$ then
$|L_{0}(\theta)|=|L_{0}(\theta_{1})|+\ldots+|L_{0}(\theta_{s_{v_{0}}})|$
and the bound follows from the inductive hypothesis.

If $\ell_{0}\in L_{0}(\theta)$ then one has
$|L_{0}(\theta)|=1 + |L_{0}(\theta_{1})|+\ldots+|L_{0}(\theta_{s_{v_{0}}})|$,
so that,  again by the inductive hypothesis,
$4|L_{0}(\theta)| \le 3|E(\theta)| - 3r_{v_{0}} - 4\,(s_{v_{0}}-1)$.
If either $r_{v_{0}}+s_{v_{0}}\ge 3$ or
$r_{v_{0}}+s_{v_{0}}=2$ and $s_{v_{0}}= 2$, the bound follows.
If $r_{v_{0}}+s_{v_{0}}=2$ and $s_{v_{0}}=1$, then $v_{0}\in V_{3}(\theta)$,
so that either $\nn_{\ell_{1}}=\vzero$ or
$2C_{0}|\oo\cdot\nn_{\ell_{1}}| \ge {\rm e}^{-\xi\mu_{v_{0}}/16}$,
by Lemma \ref{lem:3.4}, because $\ell_{0}\in L_{0}(\theta)$ and hence
$2C_{0}|\oo\cdot\nn_{\ell_{0}}|<{\rm e}^{-\xi\mu_{v_{0}}/16}$.
Therefore $\ell_{1}\notin L_{0}(\theta)$. If $v_{1}$ is the node which $\ell_{1}$
exits, call $\theta_{1}',\ldots,\theta_{s_{v_{1}}}'$ the trees whose root lines
enter $v_{1}$: one has
$|L_{0}(\theta)|=1+|L(\theta_{1}')|+\ldots+|L_{0}(\theta_{s_{v_{1}}}')|$
and hence, by the inductive hypothesis,
$4|L_{0}(\theta)| \le 3|E(\theta)|- 3r_{v_{0}}-3r_{v_{1}} -4\, (s_{v_{1}}-1)$,
where $3r_{v_{0}}+3r_{v_{1}}+4s_{v_{1}}-4 \ge 5$,
so that the bound follows in this case too.
\EP
%%%%%%%%%%%%%%%%%%%%%%%%%%%%%%%%%%%%%%%%%%%%%%%%%

%%%%%%%%%%%%%%%%%%%%%%%%%%%%%%%%%%%%%%%%%%%%%%%%%
\begin{lemma} \label{lem:3.6}
For any $k\ge 1$ and $\nn\in\ZZZ^{d}$ and any tree $\theta\in\gotT_{k,\nn}$
one has
\begin{equation} \nonumber
\left| \Val(\theta,\e) \right| \le A_{0}^{k} c_{0}^{-k} 
|\e|^{1+\frac{k+1}{8}} \prod_{v\in E(\theta)} {\rm e}^{-5\xi |\nn_{v}|/8} ,
\end{equation}
with $A_{0}$ a positive constant depending on $\Phi$, $\Gamma$ and $\rho$,
and $c_{0}$ as in Lemma \ref{lem:3.1}.
\end{lemma}
%%%%%%%%%%%%%%%%%%%%%%%%%%%%%%%%%%%%%%%%%%%%%%%%%

%%%%%%%%%%%%%%%%%%%%%%%%%%%%%%%%%%%%%%%%%%%%%%%%%
\prova
One bounds (\ref{eq:2.3}) as
\begin{equation} \nonumber
\left| \Val(\theta,\e) \right| \le |\e|^{d(\theta)}  
\Biggl( \prod_{v\in V(\theta)} \!\!\!\! |a_{p_{v}} | \Biggr) 
\Biggl( \prod_{v \in E(\theta)}  | f_{\nn_{v}} | \Biggr)
\Biggl(  \prod_{\ell \in L(\theta)} |\calG_{\ell}| \Biggr) .
\end{equation}
We deal with the propagators by using Lemma \ref{lem:3.1} as follows. 
If $\ell$ exits a node $v\in \ol{V}_{2}(\theta)$, then we have
\begin{equation} \nonumber
\left| \calG_{\ell} \right|
\!\!\!\!\! \prod_{w\in E(\theta,v)} \!\!\!\!\! 
|f_{\nn_{w}}|\,|\calG_{\ell_{w}}| \le
\frac{1}{c_{0}|\oo\cdot\nn_{\ell}|^{2}} \!\!\!
\prod_{w\in E(\theta,v)} \frac{|f_{\nn_{w}}|}{c_{0}|\oo\cdot\nn_{w}|^{2}} \le
c_{0}^{-1}C_{0}^{2}(c_{0}^{-1}C_{0}^{2}\Phi)^{|E(\theta,v)|} \!\!\!\!\!\!
\prod_{w\in E(\theta,v)} \!\!\!\!\!\! {\rm e}^{-3\xi|\nn_{w}|/4} ,
\end{equation}
where $\ell_{w}$ denotes the line exiting $w$.
For the other lines in $L_{1}(\theta)$ we distinguish three cases:
given a node $v \in V_{3}(\theta)$ and denoting by $v'$
the node which the line $\ell_{v}'$ exits, (1) if either
$\ell_{v}'\notin L_{1}(\theta,v)$ or $\ell_{v}'\in L_{1}(\theta,v')$,
we proceed as for the nodes $v\in \ol{V}_{2}(\theta)$
with $\ell=\ell_{v}$ and obtain the same bound;
(2) if $L_{1}(\theta,v)=\{\ell_{v}'\}$ and  $\ell_{v}'\notin L_{1}(\theta,v')$,
we proceed as for the nodes $v\in \ol{V}_{2}(\theta)$ with $\ell=\ell_{v}'$
and we obtain the same bound once more;
(3) if both lines $\ell_{v},\ell_{v}'$ belong to $L_{1}(\theta,v)$
and $\ell_{v}'\notin L_{1}(\theta,w)$, we bound
\begin{equation} \nonumber
\left| \calG_{\ell_{v}} \calG_{\ell_{v}'}\right|
\!\!\!\!\! \prod_{w\in E(\theta,v)} \!\!\!\!\! 
|f_{\nn_{v}}| \, |\calG_{\ell_{w}}| \le
c_{0}^{-2} C_{0}^{4}(c_{0}^{-1} C_{0}^{2}\Phi)^{|E(\theta,v)|} \!\!\!\!\!\!
\prod_{w\in E(\theta,v)} \!\!\!\!\!\! {\rm e}^{-5\xi|\nn_{w}|/8} .
\end{equation}
For all the other propagators we bound
(1) $|\calG_{\ell}| \le 1/|a|$ if $\ell$ exits a node $v\in V_{0}(\theta)$,
(2) $|\calG_{\ell}| \le c_{0}^{-1} |\oo\cdot\nn_{\ell}|^{-2}$ if $\ell$ exits an end node
and has not been already used in the bounds above for the lines
$\ell\in L_{1}(\theta)$, and
(3) $|\calG_{\ell}| \le c_{0}^{-1} |\e|^{-2}$ if $\ell \in L_{0}(\theta)$.
Then we obtain
\begin{equation} \nonumber
\left| \Val(\theta,\e) \right| \le |\e|^{d(\theta)-2|L_{0}(\theta)|}
\Gamma^{|V(\theta)|} \rho^{|N(\theta)|}
(c_{0}^{-1}C_{0}^{2})^{|V_{1}(\theta)|}
(c_{0}^{-1} C_{0}^{2}\Phi)^{|E_{1}(\theta)|} |a|^{-|V_{0}(\theta)|} 
{\rm e}^{-5\xi |\nn|/8} ,
\end{equation}
where we can bound, by using Lemma \ref{lem:3.2} and Lemma \ref{lem:3.5},
$d(\theta)-2|L_{0}(\theta)| = |E(\theta)| + |V_{1}(\theta)|-2|L_{0}(\theta)| \ge
|E(\theta)|-|L_{0}(\theta)| \ge 1 + |E(\theta)|/4\ge
1 + (k(\theta) + 1)/8$, so that the assertion  follows.
\EP
%%%%%%%%%%%%%%%%%%%%%%%%%%%%%%%%%%%%%%%%%%%%%%%%%

%%%%%%%%%%%%%%%%%%%%%%%%%%%%%%%%%%%%%%%%%%%%%%%%%
\begin{lemma} \label{lem:3.7}
For any $k\ge 1$ and $\nn\in\ZZZ^{d}$ one has
\begin{equation} \nonumber
\left| u^{(k)}_{\nn}(\e) \right| \le A_{1}^{k} c_{0}^{-k} {\rm e}^{-\xi |\nn|/2}
|\e|^{1+\frac{k+1}{8}} ,
\end{equation}
with $A_{1}$ a positive constant $C$ depending on $\Phi$, $\Gamma$, $\xi$
and $\rho$, and $c_{0}$ as in Lemma \ref{lem:3.1}.
\end{lemma}
%%%%%%%%%%%%%%%%%%%%%%%%%%%%%%%%%%%%%%%%%%%%%%%%%

%%%%%%%%%%%%%%%%%%%%%%%%%%%%%%%%%%%%%%%%%%%%%%%%%
\prova
The coefficients $u^{(k)}_{\nn}$ are given by (\ref{eq:2.4}).
Each value $\Val(\theta,\e)$ is bounded through Lemma \ref{lem:3.6}.
The sum over the Fourier labels is performed by using
a factor ${\rm e}^{-\xi|\nn_{v}|/8}$ for each end node $v\in E(\theta)$.
The sum over the other labels is easily bounded by a constant to the power $k$.
\EP
%%%%%%%%%%%%%%%%%%%%%%%%%%%%%%%%%%%%%%%%%%%%%%%%%

Lemma \ref{lem:3.7} implies that for $\e$ small enough the series (\ref{eq:2.2})
converges uniformly to a function analytic in $\pps\in\Sigma_{\xi'}$, with $\xi'<\xi/2$.
Moreover such a function is analytic in $\e \in \Omega_{\e_{0},B}$,
provided $A_{1}^{8}\e_{0}/c_{0}^{8}$ is small enough.
This completes the proof of Theorem \ref{thm:2}.

%%%%%%%%%%%%%%%%%%%%%%%%%%%%%%%%%%%%%%%%%%%%%%%%%
%%%%%%%%%%%%%%%%%%%%%%%%%%%%%%%%%%%%%%%%%%%%%%%%%
% References 
%%%%%%%%%%%%%%%%%%%%%%%%%%%%%%%%%%%%%%%%%%%%%%%%%
%%%%%%%%%%%%%%%%%%%%%%%%%%%%%%%%%%%%%%%%%%%%%%%%%

\end{document}